\documentclass[reqno]{amsart}

\usepackage{graphicx}

\newtheorem{theorem}{Theorem}[section]

\newtheorem{lemma}[theorem]{Lemma}
\newtheorem{prop}[theorem]{Proposition}

\theoremstyle{definition}
\newtheorem{definition}{Definition}[section]
\newtheorem{remark}[definition]{Remark}

\numberwithin{equation}{section}

\newcommand{\CC}{\mathbb{C}}
\newcommand{\RR}{\mathbb{R}}

\newcommand{\ZZ}{\mathbb{Z}}

\newcommand{\NN}{\mathbb{N}}
\newcommand{\AI}{A_\infty}
\newcommand{\CP}{\mathbb{C}P}
\newcommand{\RP}{\mathbb{R}P}

\newcommand{\CM}{\mathcal{M}}

\newcommand{\bx}{\boldsymbol{x}}

\begin{document}
\title[Counting real J-holomorphic discs and spheres]
{Counting real J-holomorphic discs and spheres in dimension four and six}
\author{Cheol-Hyun Cho}

\address{Department of Mathematics, Northwestern University,
Evanston, IL 60208,  Email address : cho@math.northwestern.edu }

\begin{abstract}
First, we provide another proof that the signed count of the real $J$-holomorphic 
spheres (or $J$-holomorphic discs) passing through a generic real configuration of $k$ points is independent of the choice of the real configuration and the choice of $J$, if the dimension of the Lagrangian submanifold $L$ (fixed points set of the involution) is two or three, and
also if we assume $L$ is orientable and relatively spin, and $M$ is strongly semi-positive.
This theorem was first proved by Welschinger in a more general setting, and we provide more natural approach using the degree of evaluation maps from the moduli spaces of $J$-holomorphic discs. 
Then, we define the invariant count of discs
intersecting cycles of a symplectic manifold at fixed interior marked points,
and intersecting real points at the boundary under certain assumptions.
The last result is new and was not proved by Welshinger's method.
\end{abstract}
\maketitle

\bigskip

\section{Introduction}
Let $(M,\omega)$ be a symplectic manifold with an anti-symplectic involution $\tau$
whose fixed point set is a Lagrangian submanifold $L$.
We choose a generic compatible almost complex structure $J$, which
satisfies $(\tau)_*J = -J$. We assume that $M$,$L$ and $J$ satisfies 
the above condition throughout the paper. Additionally, we assume that
$L$ is orientable and (relatively) spin, and $M$ is strongly semi-positive, which means that
for every spherical class $\alpha \in H_2(M,\ZZ)$ such that 
$[\omega](\alpha) >0$, the inequality $c_1(M)\alpha \geq 2-n$ implies
$c_1(M) \geq 1$. The first condition implies that the moduli spaces of 
$J$-holomorphic discs are orientable. The second condition
implies that the moduli spaces of simple real $J$-holomorphic (or $J_t$-holomorphic) discs form pseudo-cycles modulo disc bubbling (see section \ref{mod} for definition).

In this paper, we study the counting of 
$J$-holomorphic discs $w:(D^2,\partial D^2) \to (M,L)$ with the above assumptions 
or real $J$-holomorphic spheres $w_\CC:(\CP^1,conj) \to (M,\tau)$ preserving real structures.
We call them real $J$-holomorphic discs (or spheres) from now on.

{\bf Theorem}
{\it Assume $M,\omega,L,J$ satisfy the above assumptions.
If the dimension of $L$ is two or three, then there is a signed count of $J$-holomorphic spheres ( and discs) which passes through a given real configuration of points such that
it is independent of the choice of a real configuration of points or the choice of a generic almost complex structure. }

This result was proved by Welschinger in his beautiful but mysterious(at least to
the author) papers
\cite{W1},\cite{W2},\cite{W3} in a more general setting. The purpose of this paper is to provide conceptually more natural approach to the 
proof of the invariance of the counting using the moduli spaces of $J$-holomorphic discs. 

First of all, one naturally expect that the count of the real $J$-holomorhic spheres should
be related to the counting of real $J$-holomorphic discs since the former
can be obtained as a complex double of the latter. 
But the exact statement with orientation is rather complicated.

The most natural way to define an intersection number is as a degree of a certain map $f:M \to N$ where $M$, $N$ are oriented manifolds (for $N$ connected).
From the basic transversality theory, the signed degree of $f$, $f^{-1}(n)$ for some $n \in N$, is well-defined, and invariant under the choice of $n \in N$ if the manifold $M$ has no boundary.
In our case, we may try to consider an evaluation map $ev_k$ from the
moduli space of real $J$-holomorphic spheres to $L^k$ and consider its signed degree.
But the moduli spaces of real $J$-holomorphic spheres are not known
to be orientable (in fact, the Proposition \ref{pr:signFOOO} suggest that
it is not orientable in general). Hence we can not define the count as
a degree of the map $ev_k$ from the moduli space of real holomorphic spheres.

Our approach is to use the moduli spaces of $J$-holomorphic discs instead,
because they carry coherent orientations when the Lagrangian submanifold $L$ is orientable 
and (relatively) spin. 
But the disc moduli spaces have boundaries (actually corners) provided by the splitting of domain into two (or more) disc (or sphere) components. 
Therefore, the signed degree of $ev_k$ for each homotopy class of the discs is not well-defined.
The idea, then is to put together several moduli spaces
of holomorphic discs to make
codimension one boundaries cancel out each other, so that the sum of signed degree is well-defined. We use the conjugation $\tau$ to make sure that
there is always a matching pair for each boundary term.
But, this also does not work, as the orientation changes of the conjugations depend
on the Maslov indices of the holomorphic discs and the number of marked points.
Hence, as we put together these moduli spaces of holomorphic discs, some boundaries
cancel out, and some of them do not cancel out.

But we show that for $dim(L)=3$, if the evaluation map images of a boundary stratum
actually form codimension one boundary of the images of the non-singular stratum,
then we show that such stratum cancels out with sign. The codimension one strata
which do not have codimension one images in fact are irrelevant to the
argument to prove the well-definedness of the degree of the evaluation map.
Hence, in this case, the signed degree(which we call $A$-count) is well-defined. 

For the case $dim(L)=2$, we need to add additional sign corresponding to the cyclic
ordering of the boundary marked points, and we show that the modified counting
(which we call $B$-count) is well-defined.

The approach taken by Welschinger is 
to define a sign (spinor state) to each preimage curve, (not the
whole moduli space), and  to prove the invariance of the signed count
as the configuration of points passes though
critical points of the evaluation map $ev_k$ or a reducible curve.
The critical points of the evaluation maps are such as the
images of cuspidal curves (See\cite{W1},\cite{K}). Hence his proof is quite analytical
as he needs to analyze the change of spinor state as the configuration
crosses the critical values of the evaluation maps.

Our proof is rather combinatorial as we focus on the cancellation of signs.
We can bypass the analysis on the transversality of the evaluation maps,
because the well-definedness of the degree can be proved without such analysis.

In section 6, we show that the cancellation arguments easily generalizes to
the counting of real holomorphic discs
which intersect arbitrary pseudo-cycles of a symplectic manifold  at the
fixed interior marked points, and intersect real points at the
boundary. The count is shown to be independent of the choices involved with the above assumptions
and if the number of fixed interior marked points satisfies certain conditions.

In the last section, we compare the Welschinger's spinor state with
the coherent orientation of $J$-holomorphic discs defined
by Fukaya, Oh, Ohta and Ono in \cite{FOOO}.
But we do not know if the signed degree defined in this paper is
the same as Welschinger invariants, but we expect they are related.
They share the properties that the counting do not change
as the configuration of points pass through images of
reducible curves, or critical values of evaluation maps.

Our method does not generalize to the higher dimensional cases
in contrast to Welschinger's results \cite{W3}. It may not
neccesarily contradict his results as there may be other
combinatorial ways to achieve cancellations in the general case (as
the case of $B$-counts), but for $dim(L) \geq 4$, the counts we define in this paper, may changes as the configuration of points passes
through reducible curves in general. The argument
used for $dim(L)=2$ and $3$ does not work because there may
be too many types of reducible curves, and it is not possible
to guarantee the cancellation of all strata with
codimension one images with sign.

{\bf Acknowledgement}
We would like to thank Seongchun Kwon for reading the draft version and 
for the helpful suggestions. We would like to thank
Jean-Yves Welschinger for explaining his results to us during his visit
to Northwestern University in 2004.
\section{Orientation of the moduli space of discs and its conjugates}
\subsection{Orientation and conjugation.}
Suppose $M,\omega,J,L$ satisfy the assumption given in the introduction.
Recall that for any $J$-holomorphic map
$w:(D^2,\partial D^2) \mapsto (M,L)$, its conjugate $J$-holomorphic map is defined as
\begin{equation}\label{eq:conj}
\widetilde{w}(z) = \tau \circ w(\overline{z}), \;\; \textrm{for} \; z \in D^2
\end{equation}
We denote the homology classes as $\beta = [w] \in H_2(M,L)$, 
and $\tau_*(\beta) = [\widetilde{w}]$.
Let us first recall the following proposition from \cite{FOOO}, which 
analyzed the orientation change under the conjugation.
The space of holomorphic discs were shown to be orientable in
\cite{FOOO} (see the Proposition \ref{orientspin}). 
\begin{prop}[\cite{FOOO}Proposition 11.5, Corollary 11.9]\label{pr:signFOOO}
The map $$\tau_*: \CM_{0,k}(\beta) \to \CM_{0,k}(\tau_*(\beta))$$
given by 
$$(w,z_0,\cdots,z_{k-1}) \mapsto (\widetilde{w},\overline{z}_0,\cdots,\overline{z}_{k-1})$$
is orientation preserving if and only if
$$ \mu_L(\beta) + 2k \equiv 0 \;\;\; mod \; 4.$$
\end{prop}
\begin{remark}
The sign of $\CM_{0,k}(\beta)$ in this proposition does not incorporate the cyclic ordering
factor, which we will explain later. For the definition of the moduli spaces used 
in this paper, see section 7.
\end{remark}
Note the we always have a commuting diagram
\begin{equation}\label{commd}
\begin{array}{ccc} \CM_{0,k}(\beta) &
\stackrel{ev_{k,\beta}}{\longrightarrow} & L^k \\
 \tau_* \downarrow &    &  \parallel \\
\CM_{0,k}(\tau_* (\beta)) & \stackrel{ev_{k,\tau_* \beta}}{\longrightarrow} & L^k
\end{array}
\end{equation}
Here $ev_{k,\beta} : \CM_{0,k}(\beta) \mapsto L^k$ is the evaluation map
at the boundary marked points. Therefore, the evaluation images exactly equal to each other and the difference is the orientations of the domain moduli spaces, which can be compared using
the above theorem.
The rought idea of the proof is as follows (see \cite{FOOO} for the exact details.)
Note that the space of holomorphic discs
of homotopy class $\beta$ has real (virtual) dimension $n+ \mu(\beta)$, where
its $n$-dimensional part of the tangent space of the moduli space
is oriented by spin structure of $L$, and the remaining part of dimension $\mu(\beta)$
carry a complex orientation. Then, note that the conjugation
$\CC^{\mu/2} \to \CC^{\mu/2}$ is orientation preserving if and only if
$\mu \equiv 4$ mod 4, and the conjugation $(\partial D^2)^k \to (\partial D^2)^k$
is orientation preserving if and only if $k$ is even. 

A good example of this is the case of $(\CP^1,\RP^1)$. Denote the upper-hemisphere
disc by $U$ and the lower-hemisphere disc by $L$ (their Maslov indices are two).
In \cite{C}, we show that for the obstruction cycle ( the case $k=1$, hence
$\mu + 2k =4$), the contributions of $U$ and $L$ does not cancel out,
and when we take the boundary of Floer cohomology (the case $k=2$,
 hence $\mu + 2k =6$), the contributions cancel out. 
 
\subsection{Cyclic ordering of marked points: $A$ and $B$ counts}
Consider the moduli space of distinct $k$ marked points on $\partial D^2$, denoted
as $\CM_{0,k}^b$. Here we assume $k\geq 3$. Here the marked points can be considered as a subset of $(\partial D^2)^k \setminus \Delta$  modulo $PSL(2:\RR)$ action where 
$$\Delta : = \{(z_0,\cdots,z_{k-1}) \in (\partial D^2)^k| z_i = z_j \; \textrm{for some}
\; i \neq j \}.$$
It consists of $(k-1)!$ connected components according to the cyclic ordering of
the marked  points on $S^1$, therefore the moduli space of $J$-holomorphic discs 
$M_{0,k}(\beta)$ also has at least $(k-1)!$ connected components. 

As they are disconnected, there are two ways to orient the whole moduli spaces
$\CM_{0,k}(\beta)$. One ways is to orient the marked points as a subset of $(\partial D^2)^k$
 with counterclockwise orientation on $\partial D^2$ ( modulo $PSL(2,\RR)$).
 
 The other way is to incorporate a factor regarding the cyclic ordering of 
 the marked points.
 Namely, in the latter case, orientation may be defined by the following relation:
 For a permutation $\sigma$  of $\{1,\cdots,k\}$ (for $k \geq 3$), the map $\sigma_*:\CM_{0,k}(\beta)
 \mapsto  \CM_{0,k}(\beta)$ $$\sigma_*\big( w,(z_1,\cdots,z_k)\big) = \big(w,
 (z_{\sigma(1)},\cdots,z_{\sigma(k)})\big).$$
is orientation preserving if and only if $\sigma$ is an even permutation.

This provides two different orientations on the moduli spaces $\CM_{0,k}(\beta)$.
Therefore, if we consider the count based on each orientation, we get two different
signed counts of $J$-holomorphic discs.
We will call the count, {\bf $A$-count} (resp. {\bf $B$-count})
if the orientation does not (resp. does) incorporate cyclic ordering factor.
We remark that the cyclic ordering factor was used to in \cite{FOOO} in
the construction of $\AI$-algebra of Lagrangian submanifolds, but
it was not used in a crucial way as $\AI$ structure is constructed using
the main component.

For the purpose of the $B$-count, it is sometimes helpful to
consider each connected pieces of the moduli space $\CM_{0,k}(\beta)$ seperately.
Hence, we will denote by $\CM_{0,k}^{main}(\beta)$ the connected component where
the marked point are cyclically ordered counter clockwise as $(z_0,\cdots,z_{k-1})$,
and by $\CM_{0,k}^{\sigma(main)}(\beta)$ if the 
marked points $(z_{\sigma(0)},\cdots,z_{\sigma(k-1)})$ are cyclically ordered 
for a permutation $\sigma$ of the index set $\{ 0,\cdots,k-1\}$.

\subsection{Comparison of signs of the discs and its conjugate in the counting}
For any $J$-holomorphic disc $w$ passing through
a given configuration of real points, the conjugate disc $\widetilde{w}$ also
passes through the same configuration of points. We
compare the sign of intersections for $w$ and $\widetilde{w}$ 
as we consider counting problems.

Suppose the number of $J$-holomorphic discs of Maslov index $\mu$ passing through $k$ points on the real Lagrangian submanifold is finite (for a regular $J$ as chosen in section 7).
Namely, we will assume that 
\begin{equation}\label{eq:dim}
n + \mu(\beta) -3 + k(1-n) =0.
\end{equation}
Suppose the disc $(w,(z_0,\cdots,z_{k-1}))$ intersects the configuration
of $k$ real points $\{x_0,\cdots,x_{k-1} \}$ at the corresponding marked points.
The case of real configuration will be carried out in section \ref{sec:realcf}.

We define the orientation of the intersection to be the usual preimage orientation of
\begin{equation}\label{ta}
ev_{k,\beta}^{-1}(x_1 \times \cdots \times x_{k-1}),
\end{equation}
(See \cite{GP} Chapter 3 for the standard definition of preimage orientation.)

Now, for the moduli space $\CM_{0,k}(\tau_*\beta)$, we consider
\begin{equation}\label{tb}
ev_{k,\tau_*\beta}^{-1}(x_1 \times \cdots \times x_{k-1}),
\end{equation}

The preimages (\ref{ta}) and (\ref{tb}) are the same as unoriented sets
from the diagram (\ref{commd}). Now, there may be several discs in (\ref{ta}) but
we may compare the signs at the same time because the differences of the signs only depend
on the Maslov indices of the maps and the number of marked points, which are
fixed in this case.

For the $A$-count, Proposition \ref{pr:signFOOO} can be used to
deduce that  (\ref{ta}) and (\ref{tb}) have the same signed counts
if and only if $\mu/2 + k$ is even.

For the case of the $B$-count, the conjugation also changes the cyclic ordering of the marked points to the completely reverse order. This amounts to the new
sign contribution $(-1)^{(k-2)(k-1)/2}$.
Hence, for the $B$-count, (\ref{ta}) and (\ref{tb})
have the same signed counts if $\mu/2 + k + (k-2)(k-1)/2$ is even.
 
We provide the following table according to the above criteria,
and here we consider the (mod 4) dimension of the Lagrangian submanifold and
the (mod 4) number of point-intersection conditions.
 If (\ref{ta}) and (\ref{tb}) have the same orientations, then we denote by 1 and
 the opposite orientation by $-1$. We mark $X$  for the case which does not occur
 when we consider the condition (\ref{eq:dim}).
 
 \begin{center}
\begin{tabular}{||c||c|c|c|c||c|c|c|c||}
\hline
& \multicolumn{4}{|c||}{\bf A counts} & \multicolumn{4}{|c||}{\bf B counts} \\ \hline
$dim(L)$ &  $k = 0$ & $k = 1$ & $k=2$ & $k= 3$ &  $k = 0$ & $k = 1$ & $k=2$ & $k= 3$  \\ \hline
$0 $       &   X     &  1    & X &  -1    &    X     &  1    & X &  1 \\ \hline
$1$       &   -1     &  1  & -1 & 1    &   1     &  1  & -1 & -1   \\ \hline
$2$       &    X     &   1   & X& -1   &    X     &   1   & X& 1 \\ \hline
 $3$       &    1     &  1    & 1& 1    &    -1     &  1    & 1& -1   \\ \hline
\end{tabular}
\end{center}

We only check the first row here, and the other rows can be done similarly.
Let $dim(L) = n = 4i$ for some $i \in \NN$, then by (\ref{eq:dim}),
we have $\mu = 4lk -k +3-4l$. Note that $k$ needs to be an odd number to
have $\mu$ an even number.
Hence, for the $A$-count, we have $$\mu/2 + k \equiv (k+3)/2 \;\; mod \; 2.$$
In the case of $B$-count, 
$$ \frac{\mu}{2} + k + \frac{(k-1)(k-2)}{2} \equiv \frac{k^2+1}{2} - k \equiv
0 \;\; mod \; 2.$$

We interpret the table in the following way.
For the case the entry is $(-1)$ (for example 
$dim (L) \equiv 4, k \equiv 3 \;\; mod \; 4$), the signed count of real $J$-holomorphic
discs (or spheres) adds up to zero, if we count holomorphic discs of both homotopy classes
$\beta$ and $\tau_*\beta$.
For the case the entry is $(+1)$, the signed count of real $J$-holomorphic
discs would be twice the signed count of real $J$-holomorphic spheres.
Later, we will consider the $A$-count when $dim(L)=3$ and the $B$-count when $dim(L)=2$.
One can easily notice that in both cases, counts does not cancel out 
directly by the above table.

\section{Reducible curves and cancellations}
Let $(\Sigma,w)$ be a reducible real $J$-holomorphic disc with boundary on $L$,
whose domain $\Sigma$ consists of two disc components, denoted as $\Sigma_1,\Sigma_2$,
which intersect each other at a point. We denote the restriction
of the map $w$ on $\Sigma_i$ as $w_i$ for $i=1,2$.
We denote the homology classes as 
$$\beta_i:= [w_i(\Sigma_i)] \in H_2(M,L),$$
and we have $\beta = \beta_1+ \beta_2$. The moduli space of such
reducible curves form a codimension 1 boundary strata of the moduli
space of $J$-holomorphic discs $\CM(\beta)$. We actually consider the
case with the boundary marked points, hence the moduli space
$\CM_{0,k}(\beta)$ has codimension 1 boundary strata  such as
$$\CM_{0,k_1+1}(\beta_1) \,_{ev_i} \times_{ev_0} \CM_{0,k_2+1}(\beta_2).$$

Now, we consider the conjugate holomorphic discs of the second component 
of homotopy class $\tau_*(\beta_2)$ (See (\ref{eq:conj}) for the definition). Then, we consider the fiber product $$\CM_{0,k_1+1}(\beta_1) \,_{ev_i}\times_{ev_0} \CM_{0,k_2+1}(\tau_*(\beta_2)).$$
This can be considered as a boundary of the moduli space 
$\CM_{0,k}(\beta_1 + \tau_*(\beta_2))$.
Hence, when we consider 
$$\CM_{0,k}(\beta_1 + \beta_2) \cup \CM_{0,k}(\beta_1 + \tau_*(\beta_2)),$$
the above two boundary strata disappear by cancellation in mod 2.
Similar cancellation occurs for $\CM_{0,k}(\tau_*(\beta_1) + \beta_2)$ with
$\CM_{0,k}(\tau_*(\beta_1)+ \tau_*(\beta_2))$.  Note that
there is another choice of matching them for cancellations in mod 2, namely
compare $\CM(\beta_1+\beta_2)$ with $\CM(\tau_*(\beta_1)+\beta_2)$
and $\CM(\beta_1+\tau_*(\beta_2))$ with $\CM(\tau_*(\beta_1)+ \tau_*(\beta_2))$.

But in either way, if we consider them with orientations, they do not cancel out
in general.
We proceed more precisely. We have a signed formula as the following.
\begin{equation}\label{tta}
\CM_{0,k_1+1}(\beta_1) \,_{ev_i} \times_{ev_0} \CM_{0,k_2+1}(\beta_2)
= (-1)^{\epsilon} \partial (\CM_{0,k}( \beta_1 + \beta_2)),
\end{equation}
\begin{equation}\label{ttb}
\CM_{0,k_1+1}(\beta_1) \,_{ev_i} \times_{ev_0} \CM_{0,k_2+1}(\tau_*(\beta_2))
= (-1)^{\epsilon} \partial (\CM_{0,k}( \beta_1 + \tau_*(\beta_2)),
\end{equation}
which has the same $\epsilon$ exponent. (See \ref{signfb}).

In fact as we compare 
$\CM_{0,k}(\beta_1 + \beta_2)$ and $\CM_{0,k}(\beta_1 + \tau_*(\beta_2))$,
we need to compare connected components of different cyclic ordering
because the conjugation $\tau$ on $\beta_2$ will change the cyclic ordering
of $\tau_* \beta_2$ component for $k_2 \geq 2$. More precisely,
for (\ref{tta}), suppose the marked points $(z_0,z_1,\cdots,z_{k-1})$ are cyclically 
ordered and they split into  
$$(z_0,\cdots,z_{i-1},z_*,z_{i+k_2},\cdots,z_{k-1})\;\; \textrm{and} \;\;(z_*',z_i,z_{i+1}, \cdots,z_{i+k_2-1}),$$
where $*$ denote the marked points where we glue two discs.
Let $\sigma$ be the permutation of the index set
$\{ 0,\cdots,k-1 \}$ which maps
$$(0,\cdots,k-1) \mapsto (0,\cdots,i-1,i+k_2-1,i+k_2-2,\cdots,i,i+k_2,\cdots,k-1).$$
let $\sigma_2$ be the permutation of the index set
$\{*,i,i+1, \cdots,i+k_2-1)$ which sends
$$(*,i,i+1, \cdots,i+k_2-1) \mapsto (*,i+k_2-1,\cdots,i)$$
Then, we rewrite (\ref{tta}), (\ref{ttb}) as
\begin{equation}\label{ttta}
\CM_{0,k_1+1}^{main}(\beta_1) \,_{ev_i} \times_{ev_0} \CM_{0,k_2+1}^{main}(\beta_2)
= (-1)^{\epsilon} \partial (\CM_{0,k}^{main}( \beta_1 + \beta_2)),
\end{equation}
\begin{equation}\label{tttb}
\CM_{0,k_1+1}^{main}(\beta_1) \,_{ev_i} \times_{ev_0} \CM_{0,k_2+1}^{\sigma_2(main)}(\tau_*(\beta_2))
= (-1)^{\epsilon} \partial (\CM_{0,k}^{\sigma(main)}( \beta_1 + \tau_*(\beta_2)),
\end{equation}

We will see in the next section that if the left hand sides (LHS) 
of (\ref{ttta}) and (\ref{tttb}) cancel out with sign
then we may regard that such boundary stratum do not exist.
Namely, as we consider cobordisms, as soon as the preimage
hits a boundary stratum of a moduli space (death)
and there should be a corresponding intersection to the pair (birth) and the
cobordism will continue as if there was no boundary
in the case that the signs are opposite to each other.

It is easy to see that the cancellation of (\ref{ttta}) and (\ref{tttb})
is equivalent to comparing the orientations of 
$\CM_{0,k_2+1}^{main}(\beta_2)$ and $\CM_{0,k_2+1}^{\sigma_2(main)}(\tau_*(\beta_2))$.
 Hence for the case of $A$-counts, they do cancel out if and only if 
\begin{equation}\label{cca}
\mu(\beta_2) + 2(k_2+1) \equiv 2 \;\; mod \; 4.
\end{equation}

For the case of $B$ counts, the moduli spaces are
oriented differently according to the cyclic ordering, hence
they do cancel out if and only if 
\begin{equation}\label{ccb}
\mu(\beta_2) + 2(k_2+1) + k_2(k_2-1) \equiv 2 \;\; mod \; 4. 
\end{equation}
\begin{center}
\begin{figure}[htb!]
\includegraphics[height=3in]{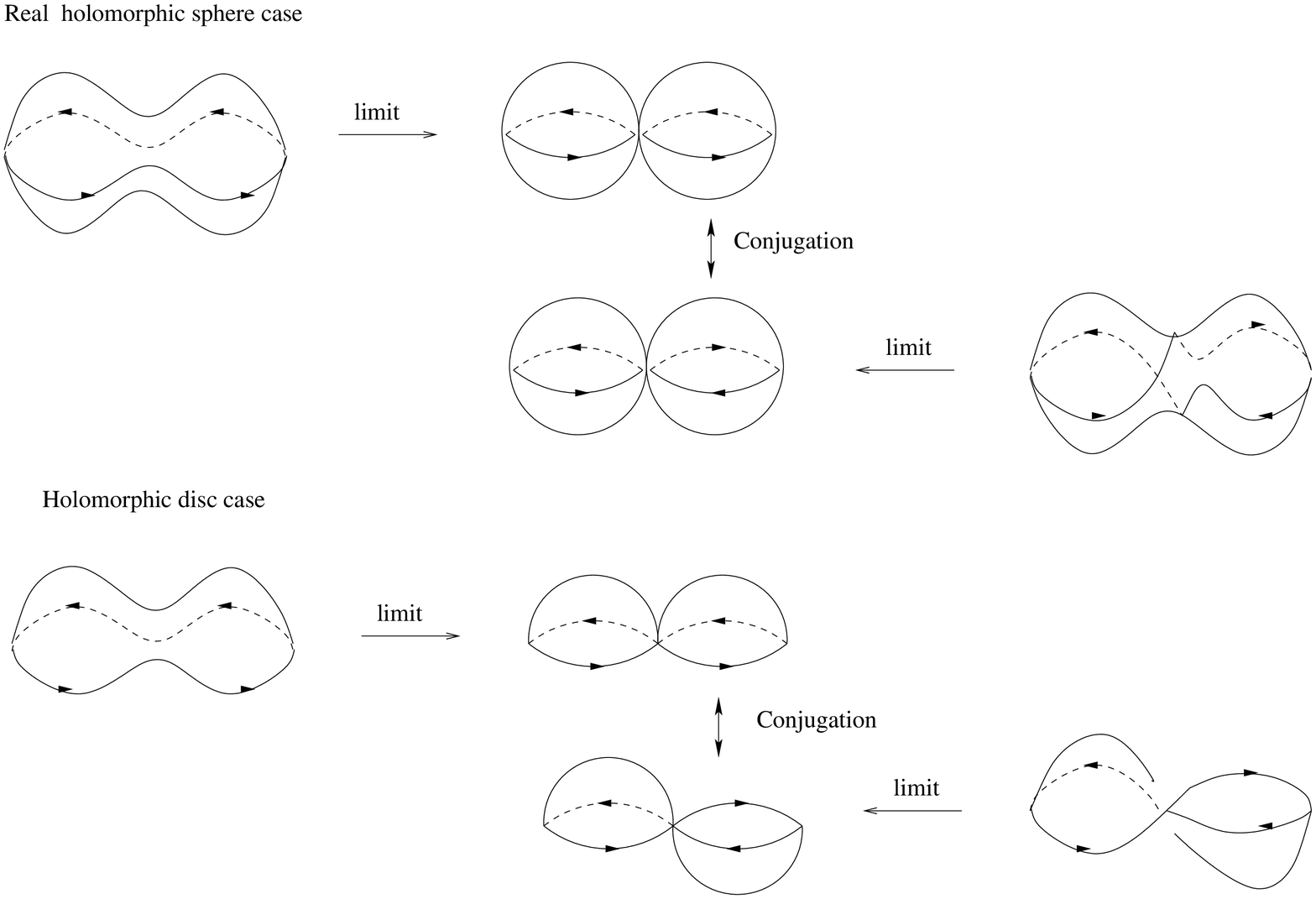}
\caption{Cancellation by conjugations}
\label{fig1}
\end{figure}
\end{center}

We remark that clearly there is an analogous story for real
$J$-holomorphic spheres as in Figure \ref{fig1}. But usually
the reducible curves are not considered as boundary of
the moduli spaces of real $J$-holomorphic curves  because of this cancellation in mod 2.
But the drawback would be that 
the moduli space of real $J$-holomorphic spheres with marked points
seems to be non-orientable in general from the above analysis.

\section{Invariance of $J$-holomorphic disc counts intersecting real points}
We define a signed count of $J$-holomorphic discs passing through generic
real $k$ points and prove its independence with respect to the choice of an almost complex structure or the choice of configuration of real points with the given assumptions.

First, we choose a spherical homology class $\alpha \in H_2(M)$, and
we say $\beta_\CC = \alpha$ if there is a holomorphic disc $w$ of class $\beta \in H_2(M,L)$
such that its complex double $w_\CC$ has a homology class $\alpha$. It is
not hard to check that if $\beta_\CC = \alpha$, then $(\tau_*\beta)_\CC = \alpha$ also.

Let $\bx = (x_1,\cdots,x_k)$ be the configuration of $k$ real points in $L$.
Now we define the count of $J$-holomorphic discs for each $\alpha \in H_2(M)$ as
$$I(\alpha,\bx) := \sharp(\bigcup_{\forall \beta, \beta_\CC = \alpha} ev_{k,\beta}^{-1}(x_1,\cdots,x_k)),$$ 
where $I=A$ or $B$ depending on the orientation of the moduli space of holomorphic discs.

\begin{theorem}\label{mainthm}
Let $M,L,J,\tau$ satisfy the assumption given in the beginning of
the paper. Then, the $A$-count when $dim(L)=3$ or the $B$-count when $dim(L)=2$
of $J$-holomorphic discs (hence real $J$-holomorphic spheres) is invariant under the choice
of a configuration of real points or the choice of a generic almost complex structure.
\end{theorem}
\begin{remark}
The case of real configuration of points is proved in section 5.
\end{remark}
\begin{proof}
To prove the independence over the choice of $k$ real points, we take a path
$$\bx(t) = (x_1(t),\cdots,x_k(t)) \in L^k$$
Consider 
\begin{equation}\label{eq:sumbeta}
\bigcup_{\forall \beta, \beta_\CC = \alpha} ev_{k,\beta}^{-1}(\bx(t)).
\end{equation}
By choosing a generic path $\bx(t)$, we may assume that $\bx(t)$ is
transversal to both the evaluation images of non-singular strata of $\CM_{0,k}(\beta)$
and the images of the codimension one strata for all homotopy classes $\beta$ with
$\beta_\CC = \alpha$.
Therefore, (\ref{eq:sumbeta}) is a collection of one dimensional manifolds with boundaries.

We will prove that the codimension one boundaries resulted from the intersection with reducible curves of (\ref{eq:sumbeta}) cancel out, and the remaining boundaries provide the cobordism between
$\bigcup_{\forall \beta, \beta_\CC = \alpha} ev_{k,\beta}^{-1}(\bx(1))$
and $\bigcup_{\forall \beta, \beta_\CC = \alpha} ev_{k,\beta}^{-1}(\bx(0))$,
hence proving the invariance of the counting.

Similarly to prove the independence over the choice of a generic almost complex structure,
we choose a generic path $\{J_t\}$ (See section 7), and consider the moduli space
$\bigcup_{0 \leq t \leq 1} \{t \} \times \CM_{0,k}(\beta,J_t)$. 
The possible boundaries arise when $t=0$ or $t=1$ or
when disc splits into reducible discs at some time $t_i \in (0,1)$.
Again, it is enough to show that the images of the boundaries of type
$ev(\{t_i\} \times \partial \CM_{0,k}(\beta,J_{t_i}))$ cancel out, to
obtain a cobordism between the case of $J_0$ and $J_1$.

First, we have the following lemma which we prove at
the end of this section.
\begin{lemma}\label{signfb}(cf. Proposition 23.2 \cite{FOOO})
$$ \partial \CM_{0,k}(\beta) = (-1)^{\epsilon} \CM_{0,k_1+1}(\beta_1)\,_{ev_i} \times _{ev_0} \CM_{0,k_2+1}(\beta_2),$$
where $\epsilon = (k_1-1)(k_2-1) + n + k_1 + (i-1)(k_2+1).$
\end{lemma}
\begin{remark}
The actual sign of this formula is not important but only 
the fact that it is the same for $\beta_2$ and $\tau_*(\beta_2)$.
Hence we postpone its proof to the end of this section. 
\end{remark}
It is also easy to note that
$$ \partial ( \{t\} \times \CM_{0,k}(\beta,J_t)) = (-1)^{\epsilon+1}
\{t\} \times \big(  \CM_{0,k_1+1}(\beta_1)\,_{ev_i} \times _{ev_0} \CM_{0,k_2+1}(\beta_2) \big).$$

There are evaluation maps at the remaining marked points
which we denote as
$$ EV_{\beta_1,\beta_2}^{k_1,k_2} : \CM_{0,k_1+1}(\beta_1)\,_{ev_i} \times _{ev_0} \CM_{0,k_2+1}(\beta_2) \to L^{k_1+k_2},$$
$$ EV_{\beta_1,\beta_2}^{k_1,k_2,t} : \{t\} \times \big( \CM_{0,k_1+1}(\beta_1,J_t)\,_{ev_i} \times _{ev_0} \CM_{0,k_2+1}(\beta_2,J_t) \big) \to L^{k_1+k_2},$$

Now, the boundary contributions resulting from the reducible curves can be written as 
\begin{equation}\label{bb}
(-1)^{\epsilon} \big(EV_{\beta_1,\beta_2}^{k_1,k_2}\big)^{-1}(\bx(t))
\end{equation}
\begin{equation}\label{bb2}
(-1)^{\epsilon +1} \big( EV_{\beta_1,\beta_2}^{k_1,k_2,t} \big) ^{-1}(\bx)
\end{equation}

It is easy to see that the difference of the orientations of (\ref{bb}) in the case of 
$(\beta_1,\beta_2)$ and the case of $(\beta_1, \tau_* \beta_2)$, is originated from the difference of orientations in $\CM_{k_2+1}(\beta_2)$ and $\CM_{k_2+1}(\tau_* \beta_2)$.
(Note that we obtain the same (\ref{bb}) in each case as unoriented sets.)
This can be done exactly as in the previous section. The same is true for (\ref{bb2}).
Hence the cancellation is determined by (\ref{cca}) and (\ref{ccb})

Therefore we only need to show that these matching pairs of
reducible curves are oriented in the  opposite way when they contribute to the cobordisms 
(\ref{bb}) or (\ref{bb2}).

Now, we first observe that to have a nontrivial preimage
either in (\ref{bb}) or (\ref{bb2}), the dimension of the moduli space $\CM_{0,k_2}(\beta_2)$ should
at least $k_2n-1$. (We restrict the evaluation map to $L^{k_2}$, and
the image should be at least codimension one to have
a non-trivial intersection in both case.)
Hence, 
\begin{equation}\label{k2eq}
 dim (\CM_{0,k_2}(\beta_2)) = n + \mu(\beta_2) - 3 + k_2 \geq k_2n -1
\end{equation}
Therefore
$$\mu(\beta_2) \geq k_2(n-1) + 2-n.$$
Then, 
\begin{eqnarray*}
\mu(\beta_1) &=& \mu(\beta) - \mu(\beta_2) \\
             & \leq & k(n-1) + 3-n - (k_2(n-1) + 2-n) \\
	     &=& k_1(n-1) + 1
\end{eqnarray*}
By applying the same argument starting from $\beta_1$,
we have obtained the following inequality once
we assume the intersection (\ref{bb}) or (\ref{bb2}) is non-trivial: 
 \begin{equation}\label{eq:ineq}
 k_i(n-1) +2-n \leq \mu(\beta_i) \leq k_i(n-1) + 1 \;\; \textrm{for each} \; i
 \end{equation}

Now we consider the cases of $dim(L)=2$ and $dim(L)=3$ seperately.
\begin{enumerate}
\item The case $dim(L) = n =2$. We show that B-count is well-defined.
(A-count is not well-defined in this case).

The above inequality becomes
$$k_i \leq \mu(\beta_i) \leq k_i +1,$$
And also we have $\mu(\beta) = k+1$ from (\ref{eq:dim}).
Note that $k$ must be odd as we assume $L$ is orientable (hence $\mu$ is even).
This implies that $k$ splits into an even number, say $k_1$ and  an odd 
number, say $k_2$. 
Then, from the above inequality, we have $\mu(\beta_1)= k_1$ and
$\mu(\beta_2) = k_2+1 $.

For the proof of the invariance, we show that the evaluation image of $\CM_{0,k_i+1}(\beta)$
has an opposite orientation as the evaluation image of $\CM_{0,k_i+1}(\tau_*\beta)$ for both 
$i=0,1$. (In fact, it is enough to show only one of them for the proof)

More precisely, the change of sign due to the proposition \ref{pr:signFOOO} plus
the sign from the cyclic ordering permutation is as follows.
For even $k_1$, 
$$k_1 + 2(k_1+1) + (k_1-1)k_1 \equiv k_1^2 + 2 \;\; mod \; 4$$
But $k_1^2$ is always a multiple of 4, hence from the proposition,
it is orientation reversing.
For odd $k_2$, 
$$k_2+1 + 2(k_2+1) + (k_2-1)k_2 \equiv k_2^2 + 1 \;\; mod \; 4$$
But $k_1^2$ is always always of type $4j+1$, hence from the proposition, it is orientation reversing. This proves the invariance of counting for $dim(L) = 2$.

\item The case $dim(L) = n =3$.
In this case we show that $A$-count is well-defined.
From the inequality (\ref{eq:ineq}), we have
$$ 2k_i-1 \leq \mu(\beta_i) \leq 2k_i+1$$
But, since $\mu(\beta_i)$ is an even number, we have
$\mu(\beta_i) = 2k_i$ for $i=1,2$.
It is enough to show that we have $\mu(\beta_i)/2 + k_i +1$ is odd.
But, it is obvious since $\mu(\beta_i)/2+k_i+1 = 2k_i +1$ is always odd.
 This finishes the proof of invariance for the case $dim(L) =3$.

\end{enumerate}
\end{proof}

\begin{proof}
This is the proof of the lemma \ref{signfb}.
We modify the result given in Proposition 23.2 \cite{FOOO} which
can be stated as
$$ \partial \CM_{0,k}(\beta) = (-1)^{(k_1-1)(k_2-1) + n + k_1}
\CM_{0,k_1+1}(\beta_1)\,_{ev_1} \times _{ev_0} \CM_{0,k_2+1}(\beta_2).$$
The only difference here from our lemma is that the $\beta_2$ disc intersect
$\beta_1$ disc at the first marked point, where as in our case
it intersects at the $i$-th marked point.
(We assume the readers are familiar to the Proposition 23.2 for this proof) 
 To remedy it,
We first reorder the orientation of the marked points of $\CM_{0,k}(\beta)$  as
$$(\partial D^2_0)\cdots(\partial D^2_{k-1}) \mapsto
(\partial D^2_{i-1})\cdots(\partial D^2_{k-1})(\partial D^2_{0})\cdots
(\partial D^2_{i-2}),$$
where $\partial D^2_i$ denotes the $i$-th marked point for $i=0,\cdots,k-1$.
  This provides a sign $(-1)^{(i-k+1)(i-1)}$
Now, we apply the proposition which now splits at $i$-th marked point of
$\beta_1$ disc. Then, the ordering of the marked points of the resulting space
$\CM_{0,k_1+1}(\beta_1)$ is 
$$(\partial D^2_{i-1})(\partial D^2_{i})(\partial D^2_{i+k_2})\cdots
(\partial D^2_{i-2}).$$
To return to the original order of marked points, we need 
additional sign 
$$(-1)^{(i-1)(2+k-i-k_2+2)}.$$
Hence by summing up the two sign factors, we prove the lemma.
\end{proof}

\section{The case of real configuration of points}\label{sec:realcf}
The main invariance theorem can be extended to the case of real configuration of points.
Let $\bx= \{ x_1,\cdots,x_k\}$ be the configuration of $k$ points on the symplectic manifold $M$.
We call $\bx$ is the real configuration of $k$ points if
$$\{\tau(x_1),\cdots,\tau(x_k)\}= \{ x_1,\cdots,x_k\}.$$
We denote by $r$ the number of real points $x_j$ with $\tau(x_j) = x_j$, and denote by $2c = k-r$.
We assume that there are at least one real point in $\bx$, namely $r \geq 1$ because
we can avoid the real holomorphic curve with empty real locus as explained \cite{W}
or to avoid (real) codimension one boundary stratum which consists of a sphere bubble attach to a disc of zero homology class as to be explained in the revision of \cite{FOOO}.

 Then, there are $c$ conjugate pairs in this real configuration of points.
Now, we assume without loss of generality that the real configuration is
arranged so that $(x_1,x_2),\cdots,(x_{2c-1},x_{2c})$ form conjugate pairs.
Then, for any such real configuration $\bx$,
we consider 
$$X(\bx) := \{(t_1,\cdots,t_c,x_{2c+1},\cdots,x_k)| t_j = x_{2j} \;\; \textrm{or} \; t_j = x_{2j-1}
= \tau(x_{2j}) \;\; \textrm{for} \; 1\leq j \leq c\}$$
We define a count of $J$-holomorphic discs which intersect $r$
real points at the boundary and intersect either of $c$ conjugate pairs at the interior.
Consider the evaluation map $Ev_{(c,r,\beta)}:\CM_{c,r}(\beta) \to  M^c \times L^r$ defined by,
$$(u,z_0^+,\cdots,z_{c-1}^+, z_0,\cdots,z_{r-1}) \mapsto  
\big( u(z_0^+), \cdots,u(z_{c-1}^+),u(z_0),\cdots,u(z_{r-1}) \big).$$
Then, we define the signed counting as
$$I(\alpha,\bx) = \sharp(\bigcup_{\forall \beta, \beta_\CC = \alpha} (Ev_{(c,r,\beta)})^{-1}(X(\bx)),$$
where $I=A$ or $B$ depending on the sign of the moduli spaces.
It is not hard to check that such a count is finite if the equation (\ref{eq:dim}) holds.

\begin{prop}
Theorem \ref{mainthm} holds also for the real configuration of points.
\end{prop}
\begin{proof}
To compare the change of an orientation under the conjugation, we consider the following commuting diagram, which is a modified version of (\ref{commd})
\begin{equation}\label{commd2}
\begin{array}{ccc} \CM_{c,r}(\beta) &
\stackrel{Ev_{(c,r,\beta)}}{\longrightarrow} &  M^c \times L^r \\
 \tau_* \downarrow &    &  ( \tau_c, id) \downarrow \;\;\;\;\;\;\;\;\;\;\;\; \\
\CM_{c,r}(\tau_*(\beta)) & \stackrel{Ev_{(c,r,\tau_*\beta)}}{\longrightarrow} &  M^c \times L^r
\end{array}
\end{equation}
Here $\tau_c:M^c \to M^c$ defined by 
$(x_1,\cdots,x_c) \mapsto \big(\tau(x_1),\cdots,\tau(x_c)\big)$ and 
$$\tau_*(w,z_0^+,\cdots,z_{c-1}^+, z_0,\cdots,z_{r-1}) = 
 ( \widetilde{w},\overline{z_0^+},\cdots,\overline{z_{c-1}^+}, \overline{z_0},\cdots,
\overline{z_{r-1}}).$$
\begin{lemma}
The map $\tau_*$ is orientation preserving if and only if
$\mu(\beta) + 2c_1 + 2r \equiv 0$  mod $4$.
And the map $\tau_c$ is orientation preserving if and only if
$c \cdot dim(M) \equiv 0$ mod $4$.
\end{lemma}
The proof of the lemma follows from the proposition \ref{pr:signFOOO}
with the fact that conjugation $\CC^c \to \CC^c$ is orientation
preserving if and only if $c$ is even.

It is easy to see that there is a bijection between the sets
$(Ev_{(c,r,\beta)})^{-1}(X(\bx))$ and $(Ev_{(c,r,\tau_*\beta)})^{-1}(X(\bx))$.
And we compare their orientations by considering the orientation changes
under both maps, $\tau_*$ and $(\tau_c,id)$ together from the diagram \ref{commd2}.

In the case of $A$-count with $dim(L)=3$, the both
inverse images carry the same orientations since
$$\mu(\beta) + 2c + 2r + 6c = \mu(\beta)+ 2k \equiv 0 \;\; mod \;4.$$
The last line follows from the table in section 2.3.
In the case of $B$-count with $dim(L)=2$, the both
inverse image carry the same orientation as
$$\mu(\beta) + 2c + 2r + 0 + (r-2)(r-1) \equiv r^2 + 3 \equiv 0 \;\; mod\; 4.$$
Here we use the fact that $\mu(\beta) = k+1$ for $dim(L)=2$.

Hence the same table as in the section 2.3 works for these cases,
which shows that the contribution from the conjugating discs do not
cancel out in these cases.

Now, we consider the proof of the invariance of the counting. 
We first show that the we obtain the same inequality (\ref{eq:ineq}) 
for the case relevant to the proof. Suppose
in the limit we have a stable disc with the domain consists of
two disc components $\Sigma_1$ and $\Sigma_2$. We assume that each component has
$c_i$ interior marked points and $r_i$ boundary marked points.
we consider $Ev_{(c_2,r_2,\beta_2)}:\CM_{c_2,r_2}(\beta) \to M^{c_2} \times L^{r_2}$.
(The case we vary an almost complex structure $J_t$ can be done analogously
and we omit the proof in that case).
Now, we observe as before that this stable limit curve will contribute
to the counting if the image of $Ev_{(c_2,r_2,\beta_2)}$ is at most of codimension one.
Hence, we have
$$n+ \mu(\beta_2) - 3 + 2c_2 + r_2 \geq c_2 \cdot 2n + r_2 \cdot n -1 = k_2n-1.$$
Note that this is the same inequality as (\ref{k2eq}), hence we obtain
the same inequality (\ref{eq:ineq}).

Now we can show the cancellation of the image of reducible curves
as in the proof of the main theorem.  We compare the sign of $Ev_{(c_i,r_i+1,\beta_i)}$
and $Ev_{(c_i,r_i+1,\tau_*\beta_i)}$ upon the condition (\ref{eq:ineq}).
First, we consider the case $dim(L)=2$. We may assume without loss of
generality that $k_1$ is even and $k_2$ is odd and $\mu(\beta_1) = k_1$
and $\mu(\beta_2) = k_2+1$. The fact that they have the 
opposite orientations follow from 
$$\mu(\beta_1) + 2c_1 + 2(r_1 +1) + r_1(r_1-1) \equiv r_1^2+2r_1+2 \equiv 2 \;\; mod \;4,$$
$$\mu(\beta_2) + 2c_2 + 2(r_2+1) + r_2(r_2-1) \equiv r_2^2 + 2r_2 +3 \equiv 2 \;\; mod \;4.$$
where the last equalities follow because $r_1$ is even and $r_2$ is odd.

For the case $dim(L)=3$, they still have the opposite orientations
because
$$\mu(\beta_i) + 2c_i + 2(r_i +1) + 2c_i \equiv 2k_i + 4c_i + 2r_i+2 \equiv 4k_i + 2 \equiv 2 
\;\; mod \;4.$$
Hence we obtain the main theorem for the real configuration of points.
\end{proof}

\section{The case of real open-closed type invaraints.}
\begin{figure}\begin{center}
\includegraphics[height=0.8in]{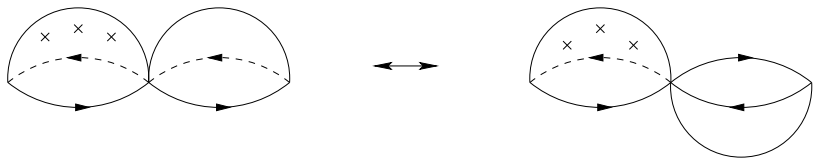}
\caption{Cancellations with fixed interior marked points}
\label{ocfig}
\end{center}
\end{figure}

In this section, we consider the count of holomorphic discs which
intersect pseudo-cycles of the symplectic manfold at the fixed interior
marked points and intersect real points at the boundary marked points.
Let $S=\{s_1,\cdots,s_l\}$ a distinct configuration of $l$ points
on the interior of the disc $D^2$, which we set to be 
fixed interior marked points. If we do not fix interior marked points,
there are too many possible strata of codimension one depending on 
partitions of the interior marked points into two discs. Because we fix the interior marked points,
all the fixed marked points stay at one disc domain in the stable limit,
so that we can apply the cancellation argument of the bubble disc component as before.
(See Figure \ref{ocfig}).

For the case $l = 1$, we obtain what we may call real open-closed Gromov-Witten invariant
as we may consider that the interior marked point is free. 
And the moduli space $\CM_{1,k}(\beta)$ also define a pseudo-cycle
modulo disc bubbling, and we will consider the obvious evaluation map 
$EV : \CM_{1,k}(\beta) \to M \times L^k$.

For the case $l \geq 2$, it may not be called Gromov-Witten type as we fix interior
marked points as in \cite{R}.
We consider evaluation maps from
the space of (parametrized) holomorphic discs $\widetilde{\CM}_{0,k}(\beta)$.
Note that we have an evaluation map
$EV : \widetilde{\CM}_{0,k}^{reg}(\beta,J) \to  M^{\times l} \times L^{\times k} $
given by $(ev_1^+,\cdots,ev_l^+,ev_0,ev_1,\cdots,ev_{k-1}),$
where $ev_j^+$ is the evaluation at the fixed interior marked point $s_j$ for $j=1,\cdots, l$
 and $ev_j$ is an evaluation at boundary marked point $z_j$. 
In general the moduli space $\widetilde{\CM}_{0,k}(\beta)$ with $EV$ map does not
define a pseudo-cycle (modulo disc bubbling).
This phenomenon was already observed in \cite{MS} in the closed case that
if we consider evaluation maps from fixed marked points, we may lose control of the positions of fixed interior marked points on
the strata containing non-simple maps. 
But for discs, it is more serious because the
structure of non-simple holomorphic disc is more complicated than that of a sphere.
This will force us to make a restrictive assumption on the number of
interior marked points related to the minimal Chern number.
This was proved in the Proposition 5.1 from \cite{C2} which is an analogue of the
theorem 5.4.1 of \cite{MS}.
\begin{prop}[Proposition 5.1 \cite{C2}]\label{prop:pseudo}
The map $EV$ defined as above, is a pseudo-cycle
modulo disc bubbling if the number of fixed interior marked points $l$
satisfies
\begin{equation}\label{cond}
l \leq min \big((\mu(\beta)-2)/2, MC \big),
\end{equation}
where $MC$ is the minimal Chern number of the manifold $M$.
\end{prop}
\begin{remark}
In the proof of the proposition in \cite{C2}, we used the fact that $J$ is integrable and
$L$ is real analytic to prove that $\mu(v_0) \geq 2MC$ in the case of (e). But in our case
of real Lagrangian submanifolds, this is satisfied without such assumptions because
$v_0$ was assumed to have more than one component in the decomposition of the
holomorphic disc, and the minimal Maslov index of $L$ is the same as $MC$ due to
the involution.
\end{remark}

Hence, we can define the counting as follows.
We choose $\vec{Q} = (Q_1,\cdots,Q_l)$, where each  $Q_i$ is a pseudo-cycle of $M$.
We define the signed count as 
$$I(\alpha,\vec{Q},\bx) = \sharp(\bigcup_{\forall \beta, \beta_\CC = \alpha} 
(Ev)^{-1}(\vec{Q},\bx)),$$
where $I=A$ or $B$ depending on the sign of the moduli spaces and when
$$ n + \mu(\beta) +k - kn - \sum_{j} deg Q_j + \delta =0.$$
Here $\delta = -1$ if $l=1$ and $\delta =0$ if $l\geq 2$.

\begin{prop}
The $A$-count, $A(\alpha,\vec{Q},\bx)$ for $dim(L)=3$ or
the $B$-count, $B(\alpha,\vec{Q},\bx)$ for $dim(L)=2$ is 
independent of the choice of the pseudo-cycles $\vec{Q}$ in their bordism classes,
the choice of the configuration of real points $\bx$, the choice of
a generic compatible complex structure( with $\tau_*J = -J$) and the choice of $S$,
if $l \leq min \big((\mu(\beta)-2)/2, MC \big)$ where $\mu(\beta)$ is the
Maslov index of the holomorphic disc to be counted.
\end{prop}
\begin{proof}
As before, this is proved using the cancellation as in Figure \ref{ocfig}.
The proof is similar to the previous cases. The analogue of
the Lemma \ref{signfb} can be proved using Lemma 7.9 of \cite{C2} to
show that orientation of the boundary strata regarding the fiber product
is the same for $\beta$ and $\tau_*(\beta)$ also in this case.
Now the inequality (\ref{eq:ineq}) for the bubble component in this case turns out to be
exactly the same as before as we can easily modify the previous proof by adding
interior intersection conditions. Hence, the cancellation with sign
can be shown as before, and this proves the proposition. In the
case with different almost complex structures and with the different choice of $S$,
the analogue of the proposition \ref{prop:pseudo} can be proved, and the result
follows from the same arguments.
\end{proof}

\section{Remark about the moduli spaces}\label{mod}
Here, we state rather well-known facts about moduli space
of real $J$-holomorphic spheres or discs.

Recall from \cite{MS} that a $k$-dimensional pseudo-cycle in $M$ is
a smooth map $f:V \to M$ defined on an
oriented $k$-dimensional manifold $V$ such that
the $dim(\Omega_f) \leq dim(V) -2$, where
$\Omega_f = \bigcap_{K \subset V,K cpt} \overline{f(V-K)}$.
Let $W$ be an oriented $k$-dimensional manifold $W$ with codimension one
boundary, then we will call a smooth map 
$F:W \to M$ a {\it pseudo-chain} if $dim(\Omega_F) \leq k-2$.
For $W$ a certain moduli space of holomorphic discs,
we say that a pseudo-chain $f:W \to M$ defines a pseudo-cycle modulo disc bubbling
if all the codimension one boundary $\partial W$ occur due to disc bubbling phenomenons on 
the moduli space $W$.

We denote by $\RR\CM(M,\alpha)^*$ the moduli space of
simple real $J$-holomorphic spheres which preserve real structures.
Namely, $u:(S^2,conj) \to (M,\tau)$ satisfies $u(conj(z)) = \tau (u(z))$.
Here $\alpha \in H_2(M)$ is a spherical homology class.

Recall that in \cite{MS}, it was shown that the moduli space of simple $J$-holomorphic
spheres form a pseudo-cycle for a generic $J$ in the strongly semi-positive case.
In the case of simple real $J$-holomorphic spheres, the similar results holds in
the following way. First, as shown in \cite{W1} the linearized $\overline{\partial}$ operator,
$D$ is $\ZZ/2\ZZ$-equivariant, where the actions are given by the involutions.
Hence for a generic $J$ with $\tau_*J = -J$, one can show that
the moduli space of simple real non-singular $J$-holomorphic spheres form
a smooth manifold by considering the equivariant verson of the proof given in \cite{MS}.
Moreover, in the case of real $J$-holomorphic curves, Gromov-compactness, and gluing of the two
real $J$-holomorphic curves work as in the case of \cite{MS}. This is because the lemma 6.1.2 of \cite{MS} can be proved $\ZZ/2\ZZ$-equivariantly as $D$ and $u$ are $\ZZ/2\ZZ$-equivariant. 
The lemma states that for the space of compatible almost complex structures $\mathcal{J}$,
the evaluation map $ev_z : \CM(A,\mathcal{J})^* \to M$ for a fixed $z \in CP^1$
is submersive, and it was crucially used to show that the moduli spaces
form a pseudo-cycle, and to make the gluing of two $J$-holomorphic curves work.

Hence, everything mentioned above also works equally well for $\CM(\beta)$, the moduli space of
$\tau$-simple $J$-holomorphic discs. Here we call $J$-holomorphic disc $w$, $\tau$-simple
if $w$ is injective away from a discrete set of points in the interior of $D^2$,
and additionally, its complex double $w_\CC$ is a simple $J$-holomorphic sphere.
($w$ is injective at $p$ if $du(p) \neq 0$ and $u^{-1}(u(p)) = \{ p \}$).
This is the moduli space which is used throughout the paper.

Now, the proof that non-simple real $J$-holomorphic spheres have images of codimension two or higher seems a bit tricky, because we do not know if there is a real version of the structure theorem of $J$-holomorphic spheres. Namely it is not clear if a non-simple real holomorphic map
can be decomposed into a simple real holomorphic map composed with a real degree $m$ map between
two $(\CP^1,conj)$'s. (Yong-Geun Oh informed us that they have a real version
of the structure theorem of $J$-holomorphic spheres in the final version of \cite{FOOO}).
But by considering the moduli space of $\tau$-simple $J$-holomorphic discs $\CM(\beta)$ and
by the structure theorem of \cite{KO} of $J$-holomorphic discs, it is not hard to
show that the images from a non $\tau$-simple real $J$-holomorphic discs is codimension two or
higher: For the case of non-simple $J$-holomorphic disc (which does not have injective points),
we can apply the structure theorem of \cite{KO} directly. For the case of a $J$-holomorphic disc $w$ which is injective away from a discrete set of points at the interior but if its complex double $w_\CC$ is not simple, one can show that the images of $w$ and $\widetilde{w}$ equal each other and in fact the image of $w$ can be considered as an image of a $J$-holomorphic sphere,
and the domain $D^2$ of $w$ is decomposed to two discs by the net $\big( D^2 \setminus w^{-1}(Im 
(w|_{\partial D^2}))\big)$ as in \cite{KO}. Hence, the images of the boundaries of such $J$-holomorphic discs can be covered by the images holomorphic discs with smaller
Maslov indices, hence their image is of codimension two or higher.

 Hence, $\CM_{0,k}(\beta)$ or $\CM_{l,k}(\beta)$ with evaluation maps
form the pseudo-cycle modulo disc bubbling. The case we vary an almost complex structure is similar and the stronly positivity 
condition implies that non-constant $J$-holomorphic discs with non-positive indices do not
appear generically.

\section{spinor state and coherent orientation of discs}\label{sec:wel}
In this section, we compare the spinor state of a real $J$-holomorphic sphere $w_\CC$
proposed by Welschinger and the coherent orientation of the $J$-holomorphic disc 
$w$ defined by Fukaya, Oh, Ohta and Ono. But we do not know whether
the Welschinger invariant is the same as the signed degree defined in this paper.

\subsection{Spinor state of a real holomorphic sphere.}
We assume that the real Lagrangian submanifold is oriented and spin
(Welschinger considers more general case when it has only pin structure). We recall the definition of spinor state defined by Welschinger in \cite{W2}, 
which was extended to general strongly positive case in \cite{W3}.
In this section, we only consider the case $n=3$.

Let $u:(\CP^1,conj) \to (M,\tau)$ be a $\ZZ/2\ZZ$-equivariant immersion.
We fix an orientation of $u|_{\RP^1}$, which is an immersed knot in $L$.
From the orientation on $\RP^1$, by choosing a disc with $\partial D^2 = \RP^1$,
we obtain a map $w:(D^2,\partial D^2) \to (M,L)$. We denote the homology class
$\beta = [w] \in H_2(M,L)$.
We have 
$$0 \to T\CP^1 \to u^*TM \to N_u \to 0.$$ 
This provides a splitting $u^*T\RR M = T\RR P^1 \oplus \RR N_u$.
The map $u$ is called {\it balanced } if its normal bundle is
a direct sum of isomorphic line bundles. i.e 
\begin{equation}\label{decom}
N_u \cong O(k-1) \oplus O(k-1).
\end{equation}

In this case, the projectivization of $N_u$ becomes 
$P(N_u) \cong CP^1 \times CP^{2}$. 
Let $x \in \RP^1$ and $(v_1,v_2)$ be the
basis of $\RR N_u|_{x}$ in the splitting (\ref{decom}).
For each $i$, there exists a unique
holomorphic line subbundle of degree $k-1$ of $N_u$ which contains $v_i$,
which projects to a constant section in $P(N_u)$. This shows
that $$\RR N_u \cong L_{\RP^1}(k-1)^{\oplus (2)}$$
If $k$ is odd, then a real line bundle $L_{\RP^1}(k-1)$
is orientable, and hence we obtain non-vanishing sections, $e_2$ and $e_3$ 
of $\RR N_u$, which are linearly independent and which comes from holomorphic sections of $N_u$. If $k$ is even, the real line bundle $L_{\RP^1}(k-1)$ is non-orientable.
But Welschinger chose a real holomorphic line bundle $G$ in $N_u$
which has degree $k-2$. (This maybe obtained using the half rotation on two
factors of $O(k-1) \oplus O(k-1)$ ). The point is that 
this $G$ provides a real holomorphic section $e_2$ of $N_u$, where its real part $\RR G$ is
orientable and $e_2|_{\RP^1}$ is non-vanishing.
If $e_2= (e_{21},e_{22})$ in the decomposition (\ref{decom}), then we set $e_3 = (-e_{22},+e_{21})$, which has the same properties as $e_2$.
Let $e_1$ be an element of the tangent space $T_w(\CM(\beta))$
which corresponds to the rotation of the domain $D^2$ at the map $u$
(holomorphic vector field along the immersed know $u|_{\RP^1}$ which is tangent
to the curve).
Then, by definition, $\{e_1,e_2,e_3\}$ form a frame of $u^*TM$ along $u|_{\RP^1}$.
We may assume that it also agree with the orientation of $L$, replacing
$e_3$ with $-e_3$ if needed.

The spinor state is defined to be $+1$ (resp. -1) if this loop
of frames lifts (resp. does not lift) to a $P_{spin(3)}(TL)$.

\subsection{Coherent orientation of discs moduli space}
We recall the coherent orientation of $J$-holomorphic discs by Fukaya-Oh-Ohta-Ono.
\begin{theorem}[\cite{FOOO}Theorem 21.1]
The moduli space of $J$-holomorphic discs is orientable, if $L
\subset(M,\omega)$ is a (relatively) spin Lagrangian submanifold.
Furthermore the choice of (relative) spin structure on $L$
determines an orientation on $\CM(L,\beta)$ canonically for all
$\beta \in \pi_2(P,L)$.
\end{theorem}

The proof is given by orienting the index bundle of the linearized
$\overline{\partial}_J$ operator. Since the zero order term does not affect
the index problem, we assume that the operator is the Dolbeault
operator $\overline{\partial}_w$. In the following proposition,
the orientation at a single disc is determined,
which can be extended to the whole moduli space with
the (relative) spin condition. (See \cite{FOOO} for details).

\begin{prop}[\cite{FOOO} Proposition 21.3]\label{orientspin}
Let $E$ be a complex vector bundle over a disc $D^2$. Let $F$ be a
totally real subbundle of $E|_{\partial D^2}$ over $\partial D^2$.
We denote by $\overline{\partial}_{(E,F)}$ the Dolbeault operator
on $D^2$ with coefficient $(E,F)$,
 $$\overline{\partial}_{(E,F)} : W^{1,p} (D^2, \partial D^2 ; E, F)
 \to L^p(D^2; E)$$
Assume F is trivial and take a trivialization of F over $\partial D^2$. Then
the trivialization gives an orientation of the virtual vector space
$Ker \, \overline{\partial}_{(E,F)} - Coker \, \overline{\partial}_{(E,F)}$
\end{prop}
\begin{proof}
Here is a proof of the proposition given in \cite{FOOO} for reader's convenience.
Suppose that the operators are surjective. (Otherwise we consider a
quotient of the target by a finite dimensional complex
subspace). 
By deforming the Hermitian
connection, we may assume that the totally real subbundle F is
trivially flat and the connection is product in a collar
neighborhood of $\partial D^2$. Let $C$ be a concentric circle in
the collar neighborhood of $\partial D^2$. If we pinch $C$ to a
point, we have the union of a disc $D^2$ and a 2-sphere $\CP^1$
with the center $O\in D^2$ and $S \in \CP^1$ identified. By the
parallel translation along radials, the trivial vector bundle F
extends up to C and its complexification gives a trivialization of
$E|_C$. Thus the bundle descends to $D^2 \cup \CP^1$. We
also denote this vector bundle by $E$. Then one can show that the
indices of the following two operators are isomorphic to each
other.
$$\overline{\partial}_{(E,F)} : W^{1,p} (D^2, \partial D^2 ; E, F)
 \to L^p(D^2; E \otimes T_{0,1}^* D^2)$$
$$\overline{\partial} : \{ (\xi_0, \xi_1) \in
W^{1,p} (D^2, \partial D^2 ; E, F) \times W^{1,p}(\CP^1,E) \,\, \vert
\, \,\xi_0(O) = \xi_1(S)\}$$
$$ \to L^p(D^2; E \otimes T_{0,1}^* D^2) \times
L^p(\CP^1; E \otimes T_{0,1}^* \CP^1)$$

We will get an orientation preserving isomorphism.
Since the real vector bundle F is trivialized, and by the above construction,
the kernel of the second operator is the kernel of the homomorphism:
\begin{equation}\label{split}
(\xi_0, \xi_1) \in  Hol (D^2,\partial D^2 : \CC^n,\RR^n)
\times Hol(\CP^1,E) \to  \xi_0(O) - \xi_1(S) \in \CC^n \cong E_S
\end{equation}
 Note that the kernel can be oriented by the orientation of
$\RR^n \cong Hol (D^2,\partial D^2 : \CC^n,\RR^n)$ since $Hol(\CP^1,E)$, and  $\CC^n$ carries a complex
orientation. This proves the Proposition. 
\end{proof}

\subsection{Relations}
As the frame $\{e_1,e_2,e_3\}$ of $u^*TL$ is made of real holomorphic sections over $\CP^1$,
we can consider them as an element of the tangent space $T_u(\CM(\beta))$.
Therefore $\{e_1,e_2,e_3\}$ provides a trivialization of $TL|_{w(\partial D^2)}$.
Note that if the homotopy class of this trivialization agrees with the
trivialization of $TL|_{w(\partial D^2)}$ obtained by the given spin structure,
then spinor state is $+1$, and otherwise $-1$. By \cite{C} (also \cite{FOOO}),
if the homotopy classes of the trivialization is the same (resp. different) than
the induced orientation on the moduli space is the same (resp. opposite).

To see exactly how this trivialization given by $\{e_1,e_2,e_3\}$
induces the pointwise orientation of the moduli space, we can
repeat the above construction in \cite{FOOO}, and
find the isomorphism between two index bundles
of $\overline{\partial}$.
From the isomorphism, we conclude that there exists a set of tangent vectors
$$\{f_1,\cdots,f_\mu \} \subset T_w(\CM(\beta))$$ carrying a canonical complex orientation such
that $\{e_1,e_2,e_3,f_1,\cdots,f_\mu \}$ is a basis of the tangent
space $T_w(\CM(\beta))$, hence gives a pointwise orientation of $T\CM(\beta)$.

Therefore Welschinger's spinor state is $+1$ (resp. -1)
if this orientation agrees (resp. disagrees) with the coherent orientation
at $T_w\CM(\beta)$ given by the spin structure.

But in order to show that Welschinger invariant is the same as the signed degree of
the evaluation map as defined in our paper, we would need to show that
$dev_k|_w:T_w\CM_{0,k}(\beta) \to T_pL^k$ is an orientation preserving map
with the orientation on $T_w\CM_{0,k}(\beta)$ given by
$\{e_1,e_2,e_3,f_1,\cdots,f_\mu \}$.
The complex orientaion on $\{f_1,\cdots,f_\mu\}$ is defined rather
abstractly, and it is not clear what would be the orientaion of
the image under evaluation map to the Lagrangian submanifold.
We do not know how to handle this problem, and we leave
it for the further research.
\bibliographystyle{amsalpha}

\end{document}